\documentclass[a4paper,11pt]{article}
\usepackage{latexsym}
\usepackage{amssymb}
\usepackage{enumerate}
\usepackage{amsfonts}
\usepackage{amsmath}
\usepackage{graphicx}
\usepackage[paper=a4paper,left=30mm,right=20mm,top=25mm,bottom=30mm]{geometry}

\newcommand{\qed}{$\Box$}

\newenvironment{@abssec}[1]{%
    \if@twocolumn

      \section*{#1}%
    \else

      \vspace{.05in}\footnotesize
      \parindent .2in
 {\upshape\bfseries #1. }\ignorespaces
    \fi}

    {\if@twocolumn\else\par\vspace{.1in}\fi}

\newenvironment{keywords}{\begin{@abssec}{\keywordsname}}{\end{@abssec}}

\newenvironment{AMS}{\begin{@abssec}{\AMSname}}{\end{@abssec}}

\newcommand\keywordsname{Key words}
\newcommand\AMSname{AMS subject classifications}
\newcommand\AMname{AMS subject classification}
\newtheorem{theorem}{Theorem}
 \newtheorem{lemma}[theorem]{Lemma}
 
 \newtheorem{proposition}[theorem]{Proposition}
\newtheorem{remark}[theorem]{Remark}
 
\def\qed{\vbox{\hrule height0.6pt\hbox{%
  \vrule height1.3ex width0.6pt\hskip0.8ex
  \vrule width0.6pt}\hrule height0.6pt
 }}

\title{A characterization of a hyperplane  in two-phase heat conductors
\thanks{This research was partially supported by the Grants-in-Aid
for Scientific Research (B) ($\sharp$ 18H01126 and $\sharp$ 17H02847)  and  JSPS Fellows ($\sharp$ 18J11430) of
Japan Society for the Promotion of Science.}}

\author{Lorenzo Cavallina\thanks{Research Center for Pure and Applied Mathematics,
Graduate School of Information Sciences, Tohoku
University, Sendai, 980-8579, Japan ({\tt cava@ims.is.tohoku.ac.jp}, {\tt  sigersak@tohoku.ac.jp}).}\  ,  
Shigeru Sakaguchi\footnotemark[2]\  \   and  
Seiichi Udagawa\thanks{Department of Mathematics, School of Medicine, Nihon University, Itabashi, Tokyo 173-0032, Japan
({\tt udagawa.seiichi@nihon-u.ac.jp}).}}

\date{}
\begin{document}
\maketitle

\begin{abstract}
We consider the Cauchy problem for the heat diffusion equation in the whole Euclidean space consisting of two media with different constant conductivities, where initially one  has temperature 0 and the other  has temperature 1.  Suppose that the interface is uniformly of class $C^6$.  We show that if  the interface  has a time-invariant constant temperature, then it must be a hyperplane. 
  \end{abstract}



\pagestyle{plain}
\thispagestyle{plain}

\section{Introduction}
\label{introduction}

Let $\Omega\subset\mathbb R^N$ be a domain with $N \ge 2$. Suppose that $\partial\Omega\not=\emptyset$ and $\partial\Omega$ is connected.
Denote by $\sigma=\sigma(x)\ (x \in \mathbb R^N)$  the conductivity distribution of the whole medium given by
\begin{equation}
\label{conductivity constants}
\sigma =
\begin{cases}
\sigma_s \quad&\mbox{in } \Omega, \\
\sigma_m \quad &\mbox{in } \mathbb R^N \setminus \Omega,
\end{cases}
\end{equation}
where $\sigma_s, \sigma_m$ are positive constants with $\sigma_s  \not= \sigma_m$. 

Let $u =u(x,t)$ be the unique bounded solution of the Cauchy problem for the heat diffusion equation:
\begin{equation}
  u_t =\mbox{ div}(\sigma \nabla u)\quad\mbox{ in }\  \mathbb R^N\times (0,+\infty) \ \mbox{ and }\ u\ ={\mathcal X}_{\mathbb R^N \setminus \Omega}\ \mbox{ on } \mathbb R^N\times
\{0\},\label{heat Cauchy}
\end{equation}
where ${\mathcal X}_{\mathbb R^N \setminus \Omega}$ denotes the characteristic function of the set $\mathbb R^N \setminus\Omega$. 

When $\partial\Omega$ is in particular a hyperplane, for instance, 
$$
\Omega = \{ x = (x_1, \dots, x_N) \in \mathbb R^N : x_1 > 0 \} \ \mbox{ and }\ \partial\Omega = \{ x = (x_1, \dots, x_N) \in \mathbb R^N : x_1 = 0 \},
$$
then we observe that 
\begin{equation}
\label{the average constant}
u(x, t) = \frac {\sqrt{\sigma_m}}{\sqrt{\sigma_s}+ \sqrt{\sigma_m}} \ \mbox{ for every } (x,t) \in \partial\Omega \times (0, +\infty).
\end{equation}
Indeed, the uniqueness of the solution of problem \eqref{heat Cauchy} yields that the solution $u$ does not depend on the variables $x_2, \dots, x_N$. The heat kernel for $N =1$ is explicitly given by \cite[p. 478]{GOO1987Tohoku}.   Denote by $G(x_1,y_1,t)$  the heat kernel written as 
\begin{eqnarray*}
G(x_1,y_1,t) &=& \left\{E_-(x_1-y_1,t) + \frac {\sqrt{\sigma_m}-\sqrt{\sigma_s}}{\sqrt{\sigma_m}+\sqrt{\sigma_s}}E_-(x_1+y_1,t)\right\}\mathcal X_{\{x_1\le 0,y_1\le 0\}}
\\
&+& \frac {2\sqrt{\sigma_m}}{\sqrt{\sigma_m}+\sqrt{\sigma_s}} E_-\!\!\left(x_1-\frac {\sqrt{\sigma_m}}{\sqrt{\sigma_s}}y_1,t\right)\mathcal X_{\{x_1\le 0,y_1 > 0\}}
\\
&+& \left\{E_+(x_1-y_1,t) + \frac {\sqrt{\sigma_s}-\sqrt{\sigma_m}}{\sqrt{\sigma_s}+\sqrt{\sigma_m}}E_+(x_1+y_1,t)\right\}\mathcal X_{\{x_1> 0,y_1> 0\}}
\\
&+& \frac {2\sqrt{\sigma_s}}{\sqrt{\sigma_m}+\sqrt{\sigma_s}} E_+\!\!\left(x_1-\frac {\sqrt{\sigma_s}}{\sqrt{\sigma_m}}y_1,t\right)\mathcal X_{\{x_1> 0,y_1 \le 0\}},
\end{eqnarray*}
where $E_\pm(z,t)$ are the Gaussian kernels with conductivities $\sigma_s, \sigma_m$ respectively on  $\mathbb R$ given by
$$
E_+(z,t) = \left(4\pi t\sigma_s \right)^{-\frac 12}\exp\!\!\left({-\frac {z^2}{4t\sigma_s }}\right)\ \mbox{ and }\ E_-(z,t) = \left(4\pi t\sigma_m \right)^{-\frac 12}\exp\!\!\left({-\frac {z^2}{4t\sigma_m}}\right)
$$
and each $\mathcal X_{\{ \cdot \}}$ denotes the characteristic function of the set $\{ \cdot \}$. Then the value of $u$ on $\partial\Omega\times(0,+\infty)$ is explicitly given by
\begin{eqnarray*}
u(0,x_2,\dots, x_N, t) &=& \int_{-\infty}^0 G(0, y_1, t)\ dy_1 
\\
&=& \int_{-\infty}^0\left\{E_-(-y_1,t) + \frac {\sqrt{\sigma_m}-\sqrt{\sigma_s}}{\sqrt{\sigma_m}+\sqrt{\sigma_s}}E_-(y_1,t)\right\} dy_1
\\
& =& \frac {\sqrt{\sigma_m}}{\sqrt{\sigma_s}+ \sqrt{\sigma_m}}.
\end{eqnarray*}

The main purpose of the present paper is to show that the converse also holds true.
%%%%%%%%%%%%%%%%%%%%%%%%%
%%%%%%%%%%%%%%%%%%%%%%%%%
%%%%%     Main theorem; Theorem 1.1 begins       %%%%%%%
%%%%%%%%%%%%%%%%%%%%%%%%%
%%%%%%%%%%%%%%%%%%%%%%%%%
\begin{theorem}
\label{th:characterization of a hyperplane} 
Let $u$ be the solution of problem \eqref{heat Cauchy}. Suppose that $\partial\Omega$ is uniformly of class $C^6$.  If there exists a constant $k$ satisfying
\begin{equation}
\label{stationary isothermic surface}
u(x,t) = k\ \mbox{ for every } (x,t) \in \partial\Omega \times (0, +\infty),
\end{equation}
then $\partial\Omega$  must be a straight line when $N=2$ and it must be a hyperplane when $N \ge 3$.
\end{theorem}

We note that if  the solution $u$ of problem \eqref{heat Cauchy} satisfies \eqref{stationary isothermic surface} for a constant  $k$,  then  $k$ must equal $\frac {\sqrt{\sigma_m}}{\sqrt{\sigma_s}+ \sqrt{\sigma_m}}$, which is the same as in \eqref{the average constant},  by Proposition \ref{modified formula on the boundary value} in section \ref{section2}.

We mention a  remark on the case where $\sigma_s = \sigma_m$.  If $\sigma_s = \sigma_m$ and $N \ge 3$, then Theorem \ref{th:characterization of a hyperplane} does not hold.  A counterexample is given in  \cite[p. 4824]{MPStransaction2006}.  Indeed, let $\mathcal H$ be a helicoid in $\mathbb R^3$.  When $\partial\Omega =\mathcal H \times \mathbb R^{N-3}$ ($\partial\Omega=\mathcal H$ for $N=3$), by the symmetry of $\mathcal H$ the solution $u$ satisfies 
\begin{equation}
\label{just half on the helicoid}
u = \frac 12\ \mbox{ on } \partial\Omega \times (0, + \infty).
\end{equation}
For convenience, we give a proof of this fact in subsection \ref{subsec: helicoid} of the Appendices.  
Moreover, when $\sigma_s = \sigma_m$, without loss of generality when $\sigma_s = \sigma_m = 1$,  by using the results of  \cite{MPStransaction2006,  Nanalysis1995} together with the explicit representation of the solution via Gaussian kernel,  we have
%%%%%%%%%%%%%%%%%%%%%%%%%
%%%%%%%%%%%%%%%%%%%%%%%%%
%%%%%     Theorem 1.2 begins       %%%%%%%
%%%%%%%%%%%%%%%%%%%%%%%%%
%%%%%%%%%%%%%%%%%%%%%%%%%
\begin{theorem}
\label{th:characterization of a hyperplane and a helicoid in the heat equation} 
Let $u$ be the unique bounded solution of the following Cauchy problem for the heat equation:
\begin{equation}
\label{heat Cauchy one phase}
 u_t =\Delta u\quad\mbox{ in }\  \mathbb R^N\times (0,+\infty) \ \mbox{ and }\ u\ ={\mathcal X}_{\mathbb R^N \setminus \Omega}\ \mbox{ on } \mathbb R^N\times
\{0\}.
\end{equation}
Suppose that  $\partial\Omega$ is of class $C^0$.
 If there exists a constant $k$ satisfying \eqref{stationary isothermic surface},
then $\partial\Omega$  must be a straight line when $N=2$,  it must be either a hyperplane  or a helicoid when $N = 3$,  and it must be a minimal hypersurface when $N \ge 4$.
\end{theorem}

The proof of Theorem \ref{th:characterization of a hyperplane} consists of two steps. In the first step,  we show that the mean curvature of $\partial\Omega$ must vanish with the aid of the barriers for the Laplace-Stieltjes transform of the solution. These barriers are constructed in \cite{CMSarXiv2018, SarXiv2019} under the assumption that $\partial\Omega$ is uniformly of class $C^6$.  Hence,  with the aid of the interior estimates  for solutions of the minimal surface equation  we notice that  $\partial\Omega$ is uniformly of class $C^\ell$ for every $\ell \in \mathbb N$.  This fact enables us to construct more precise barriers  in view of the formal WKB approximation for the Laplace-Stieltjes transform of the solution. The second step is devoted to proving that all the elementary functions of the principal curvatures of $\partial\Omega$ must vanish with the aid of the more precise barriers. Note that  we use the fact that  $\sigma_s \not=\sigma_m$ only in the second step, that is, even if $\sigma_s = \sigma_m$, we can prove that the mean curvature of $\partial\Omega$ must vanish.

\vskip 2ex
The following sections are organized as follows. In section \ref{section2}, we quote a lemma from \cite{CMSarXiv2018} and a proposition from \cite{SarXiv2019}. 
Sections \ref{section3} and \ref{section4} are devoted to the proofs of Theorems  \ref{th:characterization of a hyperplane} and  \ref{th:characterization of a hyperplane and a helicoid in the heat equation}  respectively. We also added two Appendices at the end. In subsection \ref{subsec: helicoid} we show how \eqref{just half on the helicoid} follows from the symmetry properites of the helicoid, while in subsection \ref{subsec: max pple}, we quote a maximum principle for elliptic equations with discontinuous conductivities from  \cite{SarXiv2019} and give its proof.

%%%%%%%%%%%%%%%%%%%%%%%%%%%%%%%%%%%%%%%%%%%%%%%%%%%%%%
%%%%%%%%%%%%%%%%%%%%%%%%%%%%%%%%%%%%%%%%%%%%%%%%%%%%%%
%%%%%%%%%%%%%%             Section 2 begins        %%%%%%%%%%%%%%%%%%%%%%%%
%%%%%%%%%%%%%%%%%%%%%%%%%%%%%%%%%%%%%%%%%%%%%%%%%%%%%%
%%%%%%%%%%%%%%%%%%%%%%%%%%%%%%%%%%%%%%%%%%%%%%%%%%%%%%
\setcounter{equation}{0}
\setcounter{theorem}{0}

\section{Preliminaries}
\label{section2}
Let us introduce the distance function $\delta = \delta(x)$ of $x \in \mathbb R^N$  to $\partial\Omega$ by
\begin{equation}
\label{distance function to the boundary of the domain}
\delta(x) = \mbox{ dist}(x,\partial\Omega)\ \mbox{ for }\ x \in \mathbb R^N.
\end{equation}
We quote a lemma concerning the solutions of problem \eqref{heat Cauchy} 
from \cite[Lemma 4.1]{CMSarXiv2018}, which simply comes from the maximum principle and the Gaussian bounds for the fundamental solution of $u_t=\mbox{ div}(\sigma\nabla u)$ due to Aronson \cite[Theorem 1, p. 891]{Ar1967bams}(see also \cite[p. 328]{FaS1986arma}).  Although \cite[Lemma 4.1]{CMSarXiv2018} concerns the case where $\Omega$ is bounded,  exactly the same proof is applicable  even if $\Omega$ is unbounded. 
For $\rho > 0$, we set
$$
\Omega_\rho = \{ x \in \Omega\ :\delta(x) \ge \rho \}\ \mbox{ and }\ \Omega_\rho^c =\{ x \in \mathbb R^N\setminus\Omega\ : \delta(x) \ge \rho \}.
$$
%%%%%%%%%%%%%%%%%%%%%%%%%%%%%%%%%%%%%%%%%
%%%%%%%%%%%%%%%%%%%%%%%%%%%%%%%%%%%%%%%%%
%%%%%     Lemma 2.1 begins       %%%%%%%%
%%%%%%%%%%%%%%%%%%%%%%%%%%%%%%%%%%%%%%%%%
%%%%%%%%%%%%%%%%%%%%%%%%%%%%%%%%%%%%%%%%%
\begin{lemma} 
\label{le:initial behavior and decay at infinity} 
Let $u$ be the solution of problem \eqref{heat Cauchy} with a general conductivity
$\sigma=\sigma(x)\ (x\in \mathbb R^N)$ satisfying
$$
0 < \mu \le \sigma(x) \le M\ \mbox{ for every } x \in \mathbb R^N,
$$
where $\mu, M$ are positive constants. Then the following propositions hold true:
\begin{itemize}
\item[\rm (1)]  The solution $u$ satisfies
\begin{equation}
\label{between zero and one}
0 < u < 1\  \mbox{  in } \mathbb R^N \times (0,+\infty).
\end{equation}
\item[\rm (2)] For every $\rho > 0$, there exist two positive constants $B$ and $b$ depending only on $N, \mu, M, \sigma_s, \sigma_m$ and $\rho$ such that
\begin{eqnarray}
&0 < u(x,t) < B e^{-\frac bt}\ &\mbox{ for every } (x,t) \in \Omega_\rho  \times (0,+\infty) \nonumber
\\
\mbox{ and }\quad&0 < 1-u(x,t) < B e^{-\frac bt}\ &\mbox{ for every } (x,t) \in \Omega_\rho^c \times (0,+\infty). \nonumber
\end{eqnarray}
 \end{itemize}
\end{lemma}

 Since a proposition \cite[Proposition E]{CMSarXiv2018}, where the boundary of the domain is compact,  also plays a key role in \cite{CMSarXiv2018},  
 in \cite[Proposition 2.3]{SarXiv2019} the proposition was modified in order to deal also with the case where $\partial\Omega$ is unbounded. Denote by $B_r(x)$ an open ball
in $\mathbb R^N$ with radius $r > 0$ and centered at a point $x \in \mathbb R^N$.
 
%%%%%%%%%%%%%%%%%%%%%%%%%%%%%%%%%%%%%%%%%
%%%%%%%%%%%%%%%%%%%%%%%%%%%%%%%%%%%%%%%%%
%%%%%     Proposition 2.2 begins       %%%%%%%%
%%%%%%%%%%%%%%%%%%%%%%%%%%%%%%%%%%%%%%%%%
%%%%%%%%%%%%%%%%%%%%%%%%%%%%%%%%%%%%%%%%%
  \begin{proposition}[\cite{SarXiv2019}]
\label{modified formula on the boundary value}
Let $\Omega$ be a possibly unbounded domain in $\mathbb R^N$, and let $z_0 \in \partial\Omega$.  Assume that  there exists  $\varepsilon > 0$ such that $\partial\Omega\cap B_\varepsilon(z_0)$ is of class $C^2$ and $\partial\Omega$ divides $B_\varepsilon(z_0)$ into two connected components.  Let $\sigma=\sigma(x)\ (x\in \mathbb R^N)$ be 
a general conductivity satisfying
$$
0 < \mu \le \sigma(x) \le M\ \mbox{ for every } x \in \mathbb R^N, \mbox{ and }\sigma(x) = \begin{cases} \sigma_s &\mbox{ if } x \in B_\varepsilon(z_0) \cap \Omega,\\  \sigma_m \ &\mbox{ if } x \in B_\varepsilon(z_0) \setminus \Omega,
\end{cases}
$$
where $\mu, M, \sigma_s, $ and $\sigma_m$ are positive constants. Let $u$ be the bounded solution of  problem \eqref{heat Cauchy} for this general conductivity $\sigma$.  Then, as $t \to +0$, $u$ converges to the number $\frac {\sqrt{\sigma_m}}{\sqrt{\sigma_s}+\sqrt{\sigma_m}}$ uniformly on $\partial\Omega\cap \overline{B_{\frac 12\varepsilon}(z_0)}$.
\end{proposition} 

\noindent
{\it Proof.\ } For convenience, we mention how to reduce the present case to the case where $\partial\Omega$ is bounded and of class $C^2$.  Since  $\partial\Omega\cap B_\varepsilon(z_0)$ is of class $C^2$,  we can find a bounded domain $\Omega_*$  with $C^2$ boundary $\partial\Omega_*$ satisfying
$$
\Omega\cap \overline{B_{\frac 23\varepsilon}(z_0)} \subset \Omega_*\subset \Omega\ \mbox{ and }\   \overline{B_{\frac 23\varepsilon}(z_0)}\cap\partial\Omega\subset \partial\Omega_*.
$$
Let us define the conductivity $\sigma_*=\sigma_*(x)\ (x \in \mathbb R^N)$ by
\begin{equation}
\label{conductivity constants two-phase}
\sigma_* =
\begin{cases}
\sigma_s \quad&\mbox{in } \Omega_*, \\
\sigma_m \quad &\mbox{in } \mathbb R^N \setminus \Omega_*.
\end{cases}
\end{equation}
Let $u_*=u_*(x,t)$ be the bounded solution of problem \eqref{heat Cauchy} where $\Omega$ and $\sigma$ are replaced with $\Omega_*$ and $\sigma_*$,  respectively.
 Then, by  \cite[Proposition E]{CMSarXiv2018},  as $t \to +0$, $u_*$ converges to the number $\frac {\sqrt{\sigma_m}}{\sqrt{\sigma_s}+\sqrt{\sigma_m}}$ uniformly on $\partial\Omega\cap \overline{B_{\frac 12\varepsilon}(z_0)}$. 
 
 We observe that the difference $v=u - u_*$ satisfies
\begin{eqnarray}
&&v_t =\mbox{ div}(\sigma_*\nabla v)\quad\mbox{in }\ B_{\frac 23\varepsilon}(z_0)\times (0,+\infty), \label{heat equation initial-boundary*C}
\\ 
&& |v| < 1\  \ \qquad\qquad\mbox{ in }\ \mathbb R^N \times (0,+\infty), \label{heat bounds-2C}
\\
&&v=0  \  \quad\qquad\qquad \mbox{ on } \ B_{\frac 23\varepsilon}(z_0)\times \{0\}.\label{heat initial*C}
\end{eqnarray}
Set
$$
\mathcal N = \left\{ x \in \mathbb R^N : \mbox{ dist}(x, \partial B_{\frac 23\varepsilon}(z_0)) < \frac 1{100}\varepsilon \right\} \left(= B_{\frac {203}{300}\varepsilon}(z_0)\setminus \overline{B_{\frac{197}{300}\varepsilon}(z_0)} \right).
$$
By comparing $v$ with the solutions of the Cauchy problem for the heat diffusion equation with conductivity $\sigma_*$ and initial data $\pm 2\mathcal X_{\mathcal N}$ for a short time,  with the aid of the Gaussian bounds due to Aronson \cite[Theorem 1, p. 891]{Ar1967bams}(see also \cite[p. 328]{FaS1986arma}),  we see that
there exist two positive constants $B$ and $b$ such that
\begin{equation}
\label{exponential decay of the difference near the touching point}
|v(x,t)| \le Be^{-\frac bt}\ \mbox{ for every } (x,t) \in \overline{B_{\frac 12\varepsilon}(z_0) } \times (0,\infty).
\end{equation}
Therefore,  since  $u_*$ satisfies the conclusion,  $u$ also does. \qed

%%%%%%%%%%%%%%%%%%%%%%%%%%%%%%%%%%%%%%%%%%%%%%%%%%%%%%
%%%%%%%%%%%%%%%%%%%%%%%%%%%%%%%%%%%%%%%%%%%%%%%%%%%%%%
%%%%%%%%%%%%%%             Section 3 begins        %%%%%%%%%%%%%%%%%%%%%%%%
%%%%%%%%%%%%%%%%%%%%%%%%%%%%%%%%%%%%%%%%%%%%%%%%%%%%%%
%%%%%%%%%%%%%%%%%%%%%%%%%%%%%%%%%%%%%%%%%%%%%%%%%%%%%%
\setcounter{equation}{0}
\setcounter{theorem}{0}

\section{Proof of Theorem \ref{th:characterization of a hyperplane}}
\label{section3}

First of all, Proposition \ref{modified formula on the boundary value} yields that the constant $k$ in \eqref{stationary isothermic surface} is determined by
\begin{equation}
\label{the constant on the interface}
k = \frac {\sqrt{\sigma_m}}{\sqrt{\sigma_s}+\sqrt{\sigma_m}}.
\end{equation}
Since $\partial\Omega$ is uniformly of class $C^6$, there exist two positive numbers $r$ and $K$ such that, for every point $p \in \partial\Omega$, there exist an orthogonal coordinate system $z$ and a function $\varphi \in C^6(\mathbb R^{N-1})$ such that the $z_N$ coordinate axis lies in the inward normal direction to $\partial\Omega$ at $p$, the origin is located at $p$,\   $C^6$ norm of $\varphi$ in $\mathbb R^{N-1}$  is less than $K$,\  $\varphi(0)=0,\ \nabla \varphi(0) = 0$ and the set $B_r(p) \cap\Omega$  is written as in the $z$ coordinate system
$$
 \{ z \in B_r(0) : z_N > \varphi(z_1, \dots, z_{N-1}) \}.
$$
Since $\partial\Omega$ is uniformly of class $C^6$ as explained above,  by choosing a number $\delta_0 > 0$ sufficiently small and setting  
\begin{equation}\label{inner tubular neighborhood of Omega}
\mathcal N_- = \{ x \in \Omega\ :\ 0< \delta(x) < \delta_0 \}\ \mbox{ and }\ \mathcal N_+ = \{ x \in \mathbb R^N \setminus \overline{\Omega}\ :\ 0< \delta(x) < \delta_0 \},
\end{equation}
where $\delta(x)$ is the distance function given by \eqref{distance function to the boundary of the domain}, 
we see that
\begin{eqnarray}
&&\sigma = \begin{cases} \sigma_s &\mbox{ in } \mathcal N_-,\\  \sigma_m &\mbox{ in } \mathcal N_+\end{cases}, \label{conductivity changes clearly}
\\
&& \delta \in C^6(\overline{\mathcal N_\pm}), \ \sup\left\{ \left|\frac {\partial^\alpha \delta}{\partial x^\alpha}(x)\right| : x \in \overline{\mathcal N_\pm}, |\alpha| \le 6 \right\} < +\infty,
 \label{c6 regularity}
\\
&&\mbox{ for every } x \in \overline{\mathcal N_\pm} \mbox{ there exists a unique }z = z(x) \in\partial\Omega \mbox{ with } \delta(x) = |x-z|, \qquad\label{the nearest point z from x}
\\
&& z(x) = x -\delta(x)\nabla\delta(x)\ \mbox{ for all } x \in \overline{\mathcal N_\pm}, \qquad\label{ the point z and distance from x} 
\\
&& \max_{1\le j \le N-1}|\kappa_j(z)| < \frac 1{2\delta_0}\ \mbox{ for every } z \in \partial\Omega, \qquad\label{upper bound of the curvatures on the boundary}
\end{eqnarray}
where  $\kappa_1(z), \dots, \kappa_{N-1}(z)$ denote the principal curvatures of $\partial\Omega$ at a point $z \in \partial\Omega$ with respect to the inward normal direction to $\partial\Omega$.  It is shown in \cite[Lemmas 14.16 and 14.17, p. 355]{GT1983} that 
\begin{equation}
\label{distance functions and principal curvatures}
|\nabla \delta(x)| = 1\ \mbox{ and }\ \Delta \delta(x) = \begin{cases}  - \sum\limits_{j=1}^{N-1}\frac {\kappa_j(z(x))}{1-\kappa_j(z(x)) \delta(x)} \ &\mbox{ for } x \in \mathcal N_-,
\\
   \sum\limits_{j=1}^{N-1}\frac {\kappa_j(z(x))}{1+\kappa_j(z(x)) \delta(x)} \ &\mbox{ for } x \in \mathcal N_+.
   \end{cases}
\end{equation}
We introduce elementary functions of the principal curvatures at $z \in \partial\Omega$ by
\begin{equation}
\label{elementary functions of the principal curvatures}
H_i(z) = \sum_{j_1<\cdots<j_i}\kappa_{j_1}(z)\cdots\kappa_{j_i}(z)\ \mbox{ for } i = 1, \dots, N-1,
\end{equation}
where $\frac 1{N-1}H_1(z)$ corresponds to the mean curvature of $\partial\Omega$ at $z \in \partial\Omega$ with respect to the inward normal direction to $\partial\Omega$.
Then we notice that,  for every $i =1, \dots, N-1,$ the composite function $H_i=H_i(z(x))$ satisfies that for $x \in \overline{\mathcal N_\pm}$ 
\begin{equation}
\label{differentiability of elementary functions}
H_i \in C^4(\overline{\mathcal N_\pm}), \ \sup\left\{ \left|\frac {\partial^\alpha H_i(z(x))}{\partial x^\alpha}\right| : x \in \overline{\mathcal N_\pm}, |\alpha| \le 4 \right\} < +\infty
\end{equation}
and
\begin{equation}
\label{derivative in normal direction must vanish}
\nabla\delta(x)\cdot\nabla H_i(z(x)) = 0 \ \mbox{ for  }  x \in \overline{\mathcal N_\pm}.
\end{equation}
Moreover, as in the proof of  \cite[Theorem 1.1]{SarXiv2019},  by introducing an increasing sequence of bounded subdomains in each of $\mathcal N_\pm$ together with an increasing sequence of bounded harmonic functions on each of the subdomains,  we can construct a function $\psi = \psi(x)$,  as the limit of the sequence,  on each of $\mathcal N_\pm$  satisfying 
\begin{equation}
\label{auxiliary harmonic functions}
\Delta \psi = 0\ \mbox{ in } \ \mathcal N_\pm,\ \psi = 0 \ \mbox{ on } \partial\Omega,\  \psi = 2\ \mbox{ on } \partial\mathcal N_\pm \setminus \partial\Omega  \mbox{ and }  0 < \psi < 2\ \mbox{ in } \mathcal N_\pm,
\end{equation}
even if each of $\mathcal N_\pm$ is unbounded. 
\begin{figure}[h]
\centering
\includegraphics[width=0.45\linewidth]{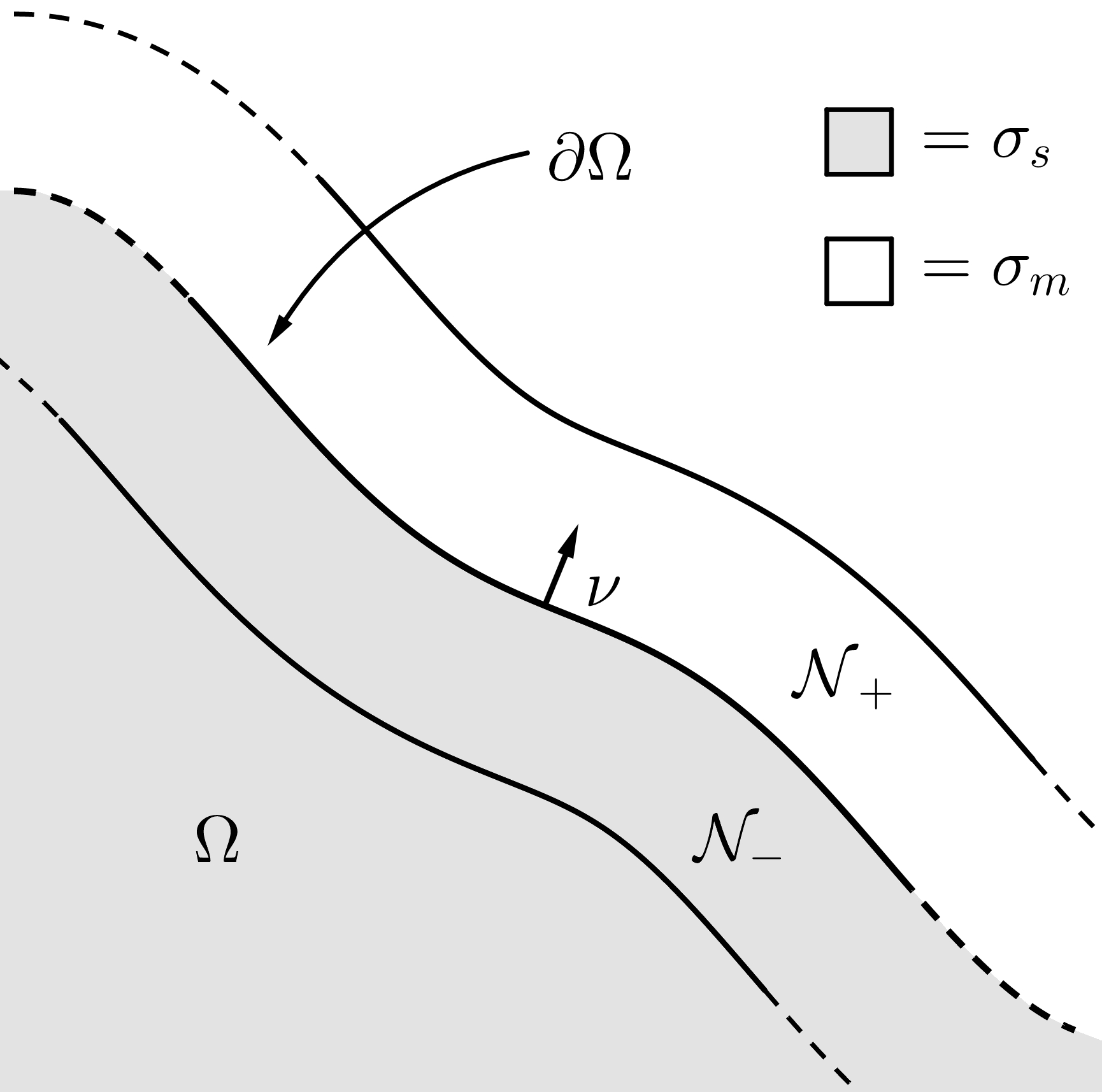}
\caption{The geometric setting used in the proof.} 
\label{Npm}
\end{figure}

As in the proofs of \cite[Theorem 1.5 in section 5]{CMSarXiv2018}, we introduce the function $w = w(x,\lambda)$ by the Laplace-Stieltjes transform of $u(x,\cdot)$  restricted on the semiaxis of real positive numbers
$$
w(x,\lambda) = \lambda\int_0^\infty\!\! e^{-\lambda t}u(x,t)\, dt \
\mbox{ for } (x,\lambda) \in \mathbb R^N\times(0,+\infty).
$$
Observe from \eqref{conductivity constants}, \eqref{heat Cauchy}, \eqref{stationary isothermic surface} and \eqref{the constant on the interface} that for every $\lambda > 0$
\begin{eqnarray}
& \sigma_s \Delta w - \lambda w = 0\ &\mbox{ in } \Omega \label{elliptic pde inside Omega},
\\
& \sigma_m \Delta (1-w) - \lambda (1-w) = 0\ &\mbox{ in } \mathbb R^N \setminus \overline{\Omega }\label{elliptic pde outside Omega},
\\
& 0 < w < 1 &\mbox{ in } \mathbb R^N, \label{value between 0 and 1}
\\
& w =  \frac {\sqrt{\sigma_m}}{\sqrt{\sigma_s}+\sqrt{\sigma_m}} \ \mbox{ and }\ \sigma_s\frac{\partial w}{\partial \nu}\big|_- =  \sigma_m\frac{\partial w}{\partial \nu}\big|_+&\mbox{ on } \partial\Omega, \label{interface condition for lambda}
\end{eqnarray}
where $\nu$ denotes the outward unit normal vector to $\partial\Omega$,  $+$ denotes the limit from outside of $\Omega$ and $-$ that from inside of $\Omega$. 
Moreover, it follows from (2)  of Lemma \ref{le:initial behavior and decay at infinity} that there exist two positive constants $\tilde B$ and $\tilde b$ satisfying:
\begin{eqnarray}
&0 < w(x,\lambda) \le \tilde B e^{-\tilde b \sqrt{\lambda}}\  &\mbox{ for every } (x,\lambda) \in \left(\partial \mathcal N_-\setminus\partial\Omega\right) \times (0, +\infty),\label{exp decay inside Omega}
\\
& 0 < 1-w(x,\lambda) \le \tilde B e^{-\tilde b \sqrt{\lambda}}\  &\mbox{ for every } (x,\lambda) \in \left(\partial \mathcal N_+\setminus\partial\Omega\right) \times (0, +\infty).\label{exp decay outside Omega}
\end{eqnarray}

%%%%%%%%%%%%%%%%%%%%%%%%%%%%%%%%%%%%%%%%%%%
%%%%%%%%%%%%%%%%%%%%%%%%%%%%%%%%%%%%%%%%%
%%%%%%%%%%%%%%%%%%%%%%%%%%%%%%%%%%%%%%%%%
%%%%%     Subsection 3.1 begins       %%%%%%%%
%%%%%%%%%%%%%%%%%%%%%%%%%%%%%%%%%%%%%%%%%
%%%%%%%%%%%%%%%%%%%%%%%%%%%%%%%%%%%%%%%%%

\subsection{Proving that the mean curvature of $\partial\Omega$ vanishes}
\label{subsec:vanishing mean curvature}

Let us first consider $w$ on $\mathcal N_-$.
Since $w$ satisfies \eqref{elliptic pde inside Omega} and the first equality of \eqref{interface condition for lambda},   in view of the formal WKB approximation of $w$ for  sufficiently large $\tau = \frac \lambda{\sigma_s}$
$$
w(x,\lambda) \sim \frac {\sqrt{\sigma_m}}{\sqrt{\sigma_s}+\sqrt{\sigma_m}}e^{-\sqrt{\tau}\delta(x)} \sum_{j=0}^\infty A_j(x) \tau^{-\frac j2}\ \mbox{ with some coefficients } \{ A_j(x) \}, 
$$ 
we introduce two functions $f_{1, \pm} = f_{1, \pm}(x,\lambda)$ defined for $(x,\lambda) \in \overline{\mathcal N_-} \times (0,+\infty)$ by
$$
f_{1,\pm}(x,\lambda) = \frac {\sqrt{\sigma_m}}{\sqrt{\sigma_s}+\sqrt{\sigma_m}}e^{-\frac{\sqrt{\lambda}}{\sqrt{\sigma_s}} \delta(x)}\left[A_0(x) + \frac {\sqrt{\sigma_s}}{\sqrt{\lambda}}A_{1,\pm}(x)\right],
$$
where 
\begin{eqnarray}
&&A_0(x) = \left\{\prod\limits_{j=1}^{N-1}\Bigl[1-\kappa_j(z(x))\delta(x)\Bigr]\right\}^{-\frac 12}, \label{A0}
\\
&&A_{1,\pm}(x) = \int_0^{\delta(x)}\left[\frac 12\,\Delta A_0(x(\tau)) \pm 1\right]\exp\left(-\frac 12\,\int_\tau^{\delta(x)} \Delta \delta(x(\tau')) d\tau'\right)d\tau, \nonumber
\end{eqnarray}
with $x(\tau) = z(x) - \tau\,\nu(z(x))$ for $0<\tau<\delta(x)$.   We observe that for $x \in \overline{\mathcal N_-}$
\begin{equation}
\label{expansion with higher order curvatures}
\prod\limits_{j=1}^{N-1}\Bigl[1-\kappa_j(z(x))\delta(x)\Bigr] = 1 + \sum_{i=1}^{N-1} (-1)^i H_i(z(x)) (\delta(x))^i.
\end{equation}
With \eqref{distance functions and principal curvatures},  \eqref{derivative in normal direction must vanish} and \eqref{expansion with higher order curvatures}  at hand, by straightforward computations we obtain that
\begin{equation}
\label{gradient of two functions}
\nabla\delta\cdot\nabla A_0 = -\frac 12(\Delta\delta)A_0, \quad \nabla\delta\cdot\nabla A_{1,\pm} = -\frac 12(\Delta\delta)A_{1, \pm} + \frac 12 \Delta A_0 \pm 1 \ \mbox{ in }  \ \overline{\mathcal N_-}, 
\end{equation}
\begin{equation}
\label{for super and subsolutions}
\sigma_s\Delta f_{1,\pm} - \lambda f_{1,\pm} = \frac {\sigma_s\sqrt{\sigma_m}}{\sqrt{\sigma_s}+\sqrt{\sigma_m}} e^{-\frac{\sqrt{\lambda}}{\sqrt{\sigma_s}} \delta(x)}\left(\mp 2 + \frac {\sqrt{\sigma_s}}{\sqrt{\lambda}}\Delta A_{1,\pm}\right)  \ \mbox{ in }  \ \overline{\mathcal N_-},
\end{equation}
and
\begin{equation}
\label{inner boundary Dirichlet}
A_0 = 1,\ A_{1,\pm} = 0,  \quad f_{1,\pm} = \frac {\sqrt{\sigma_m}}{\sqrt{\sigma_s}+\sqrt{\sigma_m}}  \ \mbox{ on } \ \partial\Omega, 
\end{equation}
for every $\lambda>0$.  Moreover, \eqref{c6 regularity}, \eqref{upper bound of the curvatures on the boundary},  \eqref{differentiability of elementary functions} and \eqref{expansion with higher order curvatures}  yield that 
\begin{equation}
\label{boundedness coming from c6}
|\Delta A_{1,\pm}| \le c_1\ \mbox{ in } \overline{\mathcal N_-}
\end{equation}
for some positive constant $c_1$. 
Therefore, it follows from \eqref{for super and subsolutions},  \eqref{boundedness coming from c6}, \eqref{exp decay inside Omega} and the definition of $f_{1,\pm}$ 
that there exist two positive constants $\lambda_1$ and $\eta_1$ such that 
\begin{eqnarray}
&&\sigma_s\Delta f_{1, +} - \lambda f_{1, +} < 0 < \sigma_s\Delta f_{1, -} - \lambda f_{1, -}\ \mbox{ in } \ \overline{\mathcal N_-},\label{key differential inequalities}
\\
&& \max\{ |f_{1, +}|, |f_{1, -}|, w \} \le e^{-\eta_1\sqrt{\lambda}}\ \mbox{ on } \ \partial \mathcal N_- \setminus\partial\Omega,\label{key inner boundary decay estimates}
\end{eqnarray}
for every  $\lambda \ge \lambda_1$. 

 For every $(x,\lambda) \in \overline{\mathcal N_-} \times (0,+\infty)$,  we define the two functions $w_{1, \pm} = w_{1, \pm}(x,\lambda)$  by
\begin{equation}
\label{upper and lower barriers}
w_{1, \pm}(x,\lambda) = f_{1, \pm}(x,\lambda) \pm \psi(x)e^{-\eta_1\sqrt{\lambda}},
\end{equation}
where $\psi(x)$ is given by \eqref{auxiliary harmonic functions}.
Then, in view of \eqref{elliptic pde inside Omega}, \eqref{interface condition for lambda}, \eqref{inner boundary Dirichlet}, \eqref{key differential inequalities} and \eqref{key inner boundary decay estimates}, we notice that 
\begin{eqnarray}
& \sigma_s\Delta w_{1, +} - \lambda w_{1, +} < 0 = \sigma_s\Delta w - \lambda w < \sigma_s\Delta w_{1, -} - \lambda w_{1, -}\ &\mbox{ in } \mathcal N_-,\nonumber
\\
& w_{1, +} = w =w_{1, -} = \frac {\sqrt{\sigma_m}}{\sqrt{\sigma_s}+\sqrt{\sigma_m}}\ &\mbox{ on } \partial\Omega,\label{key equality for gradient estimates}
\\
& w_{1, -} < w < w_{1, +}\ &\mbox{ on } \partial \mathcal N_- \setminus\partial\Omega,\nonumber
\end{eqnarray}
for every  $\lambda \ge \lambda_1$,  and hence we get that
\begin{equation}
\label{barriers work well}
w_{1, -} < w < w_{1, +}\ \mbox{ in } \ \mathcal N_-,
\end{equation}
for every  $\lambda \ge \lambda_1$, by the comparison principle (see Proposition \ref{prop:maximum principle on unbounded domains} in Appendix).
Thus, combining \eqref{barriers work well} with \eqref{key equality for gradient estimates}  yields that 
\begin{equation}
\label{key inequalities on the boundary}
\frac {\partial w_{1, +}}{\partial \nu}\le \frac {\partial w}{\partial\nu}\Big|_- \le \frac{\partial w_{1, -}}{\partial\nu} \ \mbox{ on } \ \partial\Omega,
\end{equation}
for every  $\lambda \ge \lambda_1$.

Therefore, by recalling the definition of $w_{1,\pm}$,  it follows from \eqref{gradient of two functions}, \eqref{inner boundary Dirichlet} and \eqref{distance functions and principal curvatures}  that,  for every  $\lambda \ge \lambda_1$,  we have the following chain of inequalities on $\partial\Omega$ :
\begin{eqnarray}
&&\frac {\sigma_s\sqrt{\sigma_m}}{\sqrt{\sigma_s}+\sqrt{\sigma_m}}\left\{-\frac 12 \sum\limits_{j=1}^{N-1}\kappa_j - \frac {\sqrt{\sigma_s}}{\sqrt{\lambda}}\left(\frac 12\Delta A_0+1\right)\right\} + \sigma_s\frac {\partial\psi}{\partial\nu} e^{-\eta_1\sqrt{\lambda}} \nonumber
\\
&&\le  \sigma_s\frac{\partial w}{\partial\nu}\Big|_-\!\! - \frac {\sqrt{\sigma_s}\sqrt{\sigma_m}}{\sqrt{\sigma_s}+\sqrt{\sigma_m}} \sqrt{\lambda}
  \nonumber
\\
&&\le \frac {\sigma_s\sqrt{\sigma_m}}{\sqrt{\sigma_s}+\sqrt{\sigma_m}}\left\{-\frac 12\sum\limits_{j=1}^{N-1}\kappa_j - \frac {\sqrt{\sigma_s}}{\sqrt{\lambda}}\left(\frac 12\Delta A_0-1\right)\right\} - \sigma_s \frac {\partial\psi}{\partial\nu} e^{-\eta_1\sqrt{\lambda}}.\qquad\label{bounds-for-distance-and-curvatures}
\end{eqnarray}
This implies that on $\partial\Omega$
\begin{equation}
\label{estimate from inside}
-\frac {\sigma_s\sqrt{\sigma_m}}{2(\sqrt{\sigma_s}+\sqrt{\sigma_m})} \sum\limits_{j=1}^{N-1}\kappa_j
= \sigma_s\frac{\partial w}{\partial\nu}\Big|_-\!\! - \frac {\sqrt{\sigma_s}\sqrt{\sigma_m}}{\sqrt{\sigma_s}+\sqrt{\sigma_m}} \sqrt{\lambda} + O(1/\sqrt{\lambda})\ \mbox{ as } \lambda \to +\infty. 
\end{equation}

Next, we consider $1-w$ on $\mathcal N_+$. By the similar arguments as above, since
 $$
 1-w = \frac {\sqrt{\sigma_s}}{\sqrt{\sigma_s}+\sqrt{\sigma_m}}\ \mbox{ on } \partial\Omega,
 $$
 we can construct barriers for $1-w$ on $\mathcal N_+$ with the aid of \eqref{exp decay outside Omega}  by replacing $\sigma_s$ with $\sigma_m$. Thus,  proceeding similarly yields that on $\partial\Omega$
 \begin{equation}
\label{estimate from outside}
\frac {\sigma_m\sqrt{\sigma_s}}{2(\sqrt{\sigma_s}+\sqrt{\sigma_m})} \sum\limits_{j=1}^{N-1}\kappa_j
= \sigma_m\frac{\partial w}{\partial\nu}\Big|_+\!\! - \frac {\sqrt{\sigma_s}\sqrt{\sigma_m}}{\sqrt{\sigma_s}+\sqrt{\sigma_m}} \sqrt{\lambda} + O(1/\sqrt{\lambda})\ \mbox{ as } \lambda \to +\infty,
\end{equation}
where we have taken into account both the sign of the mean curvature with \eqref{distance functions and principal curvatures} and the normal direction to $\partial\Omega$. 
 Therefore, by combining  \eqref{estimate from inside} and  \eqref{estimate from outside}  with the second equality of \eqref{interface condition for lambda} we conclude that on $\partial\Omega$
 $$
H_1= \sum\limits_{j=1}^{N-1}\kappa_j
= O(1/\sqrt{\lambda})\ \mbox{ as } \lambda \to +\infty, 
$$
 and hence the mean curvature of $\partial\Omega$ must vanish, that is, $\partial\Omega$ is a minimal hypersurface properly embedded in $\mathbb R^N$ (see \eqref{elementary functions of the principal curvatures} for $H_1$). In particular when $N=2$, the curvature of the curve $\partial\Omega$ vanishes and the conclusion of Theorem \ref{th:characterization of a hyperplane} holds.
 
  Note that in this subsection \ref{subsec:vanishing mean curvature} we did not use the fact that $\sigma_s \not=\sigma_m$.
 
 %%%%%%%%%%%%%%%%%%%%%%%%%%%%%%%%%%%%%%%%%%%
%%%%%%%%%%%%%%%%%%%%%%%%%%%%%%%%%%%%%%%%%
%%%%%%%%%%%%%%%%%%%%%%%%%%%%%%%%%%%%%%%%%
%%%%%     Subsection 3.2  begins       %%%%%%%%
%%%%%%%%%%%%%%%%%%%%%%%%%%%%%%%%%%%%%%%%%
%%%%%%%%%%%%%%%%%%%%%%%%%%%%%%%%%%%%%%%%%
 \subsection{Proving that all the principal curvaures of $\partial\Omega$ vanish and $\partial\Omega$ must be a hyperplane}
\label{subsec:vanishing principal curvatures}

We may consider the case where $N \ge 3$.
It suffices to show that $H_i = 0$ on $\partial\Omega$ for every $i=1, \dots, N-1$.  Since we already know in subsection \ref{subsec:vanishing mean curvature} that $H_1 = 0$ on $\partial\Omega$, we start induction with supposing that there exists a number $p \in \{ 2, \dots, N-1 \}$ satisfying
\begin{equation}
\label{assumption for induction}
H_1 = \cdots = H_{p-1} = 0\ \mbox{ on } \partial\Omega.
\end{equation}
Then we will  prove that $H_p =0$ on $\partial\Omega$.  By subsection \ref{subsec:vanishing mean curvature},  $\partial\Omega$ must be real analytic and
  moreover, by the interior estimates for solutions of the minimal surface equation (see \cite[Corollary 16.7, p. 407]{GT1983}), we see that $\partial\Omega$ is uniformly of class $C^\ell$ for every $\ell \in \mathbb N$, and hence \eqref{c6 regularity} and \eqref{differentiability of elementary functions} are improved as follows: For every $\ell \in \mathbb N$, 
 \begin{equation}
 \label{uniform c infty regularity}
 \sup\left\{ \left|\frac {\partial^\alpha \delta}{\partial x^\alpha}(x)\right| : x \in \overline{\mathcal N_\pm}, |\alpha| \le \ell \right\} < +\infty,
 \end{equation}
 and 
 \begin{equation}
\label{differentiability of elementary functions c infty}
\sup\left\{ \left|\frac {\partial^\alpha H_i(z(x))}{\partial x^\alpha}\right| : 1 \le i \le N-1,\  x \in \overline{\mathcal N_\pm}, |\alpha| \le \ell \right\} < +\infty.
\end{equation}
Therefore we can introduce the following more precise barriers $f_{n,\pm} = f_{n,\pm}(x,\lambda)$ for $w$ on $\mathcal N_-$ such that   for $(x,\lambda) \in \overline{\mathcal N_-} \times (0,+\infty)$ and for every $n \ge 2$ 
$$
f_{n,\pm}(x,\lambda) = \frac {\sqrt{\sigma_m}}{\sqrt{\sigma_s}+\sqrt{\sigma_m}}e^{-\frac{\sqrt{\lambda}}{\sqrt{\sigma_s}} \delta(x)}\left[A_0(x) + \sum_{j=1}^{n-1}\left(\frac {\sqrt{\sigma_s}}{\sqrt{\lambda}}\right)^j\!\!A_j(x) + \left(\frac {\sqrt{\sigma_s}}{\sqrt{\lambda}}\right)^n\!\!A_{n,\pm}(x)\right],
$$
where $A_0$ is given by \eqref{A0} and for  $j=1, \cdots, n-1$,
\begin{eqnarray}
&&A_{j}(x) = \int_0^{\delta(x)}\left[\frac 12\,\Delta A_{j-1}(x(\tau))\right]\exp\left(-\frac 12\,\int_\tau^{\delta(x)} \Delta \delta(x(\tau')) d\tau'\right)d\tau,\label{Aj}
\\
&&A_{n,\pm}(x) = \int_0^{\delta(x)}\left[\frac 12\,\Delta A_{n-1}(x(\tau)) \pm 1\right]\exp\left(-\frac 12\,\int_\tau^{\delta(x)} \Delta \delta(x(\tau')) d\tau'\right)d\tau\nonumber
\end{eqnarray}
with $x(\tau) = z(x) - \tau\,\nu(z(x))$ for $0<\tau<\delta(x)$.  

With \eqref{distance functions and principal curvatures}, \eqref{derivative in normal direction must vanish} and \eqref{expansion with higher order curvatures} at hand, by straightforward computations we obtain that, in $\overline{\mathcal N_-}$ (compare with \eqref{gradient of two functions}--\eqref{boundedness coming from c6}):
\begin{eqnarray}
&&\nabla\delta\cdot\nabla A_0 = -\frac 12(\Delta\delta)A_0,\label{gradient of coefficient 1}
\\
&&\nabla\delta\cdot\nabla A_{j} = -\frac 12(\Delta\delta)A_j + \frac 12 \Delta A_{j-1} \ \mbox{ for } j = 1, \dots, n-1, \label{gradient of coefficient j}
\\
&&\nabla\delta\cdot\nabla A_{n,\pm} = -\frac 12(\Delta\delta)A_{n,\pm} + \frac 12 \Delta A_{n-1} \pm 1,\label{gradient of coefficient n}
\end{eqnarray}
\begin{equation}
\label{for super and subsolutions more precise}
\sigma_s\Delta f_{n,\pm} - \lambda f_{n,\pm} = \frac {\sigma_s\sqrt{\sigma_m}}{\sqrt{\sigma_s}+\sqrt{\sigma_m}} \left(\frac {\sqrt{\sigma_s}}{\sqrt{\lambda}}\right)^{n-1}\!\! e^{-\frac{\sqrt{\lambda}}{\sqrt{\sigma_s}} \delta(x)}\left(\mp 2 + \frac {\sqrt{\sigma_s}}{\sqrt{\lambda}}\Delta A_{n,\pm}\right),
\end{equation}
and on $\partial\Omega$
\begin{equation}
\label{inner boundary Dirichlet more precise}
A_0 = 1,\  A_1=\cdots=A_{n-1}=A_{n,\pm} = 0,  \quad f_{n,\pm} = \frac {\sqrt{\sigma_m}}{\sqrt{\sigma_s}+\sqrt{\sigma_m}}, 
\end{equation}
for every $\lambda>0$.  Moreover, \eqref{uniform c infty regularity}, \eqref{differentiability of elementary functions c infty}, \eqref{upper bound of the curvatures on the boundary} and \eqref{expansion with higher order curvatures} yield that 
\begin{equation}
\label{boundedness coming from all natural numbers}
|\Delta A_{n,\pm}| \le c_n\ \mbox{ in } \overline{\mathcal N_-}
\end{equation}
for some positive constant $c_n$.  Then,  by replacing $f_{1,\pm}$ with $f_{n,\pm}$, we can use the same comparison arguments as in \eqref{key differential inequalities} - \eqref{key inequalities on the boundary} of subsection \ref{subsec:vanishing mean curvature} 
to conclude that there exist two positive constants $\lambda_n$ and $\eta_n$ satisfying
\begin{equation}
\label{key inequalities on the boundary more precise}
\frac {\partial w_{n, +}}{\partial \nu}\le \frac {\partial w}{\partial\nu}\Big|_- \le \frac{\partial w_{n, -}}{\partial\nu} \ \mbox{ on } \ \partial\Omega
\end{equation}
for every  $\lambda \ge \lambda_n$, where 
\begin{equation}
\label{upper and lower barriers more precise}
w_{n, \pm}(x,\lambda) = f_{n, \pm}(x,\lambda) \pm \psi(x)e^{-\eta_n\sqrt{\lambda}} 
\end{equation}
with $\psi(x)$ given by \eqref{auxiliary harmonic functions}. Since $\Delta \delta = 0$ on $\partial\Omega$, it follows from \eqref{distance functions and principal curvatures},  \eqref{gradient of coefficient 1}--\eqref{gradient of coefficient n} and \eqref{inner boundary Dirichlet more precise} that on $\partial\Omega$
\begin{eqnarray}
\frac {\partial w_{n, \pm}}{\partial \nu} &=& -\nabla \delta\cdot\nabla w_{n, \pm}\nonumber
\\
&=& \frac{\sqrt{\sigma_m}}{\sqrt{\sigma_s} + \sqrt{\sigma_m}} \left\{ \frac {\sqrt{\lambda}}{\sqrt{\sigma_s}} - \frac 12 \sum_{j = 1}^{n-1}\left(\frac {\sqrt{\sigma_s}}{\sqrt{\lambda}}\right)^j\Delta A_{j-1} - \frac 12 \left(\frac {\sqrt{\sigma_s}}{\sqrt{\lambda}}\right)^n (\Delta A_{n, \pm} \pm 2)\right\} \nonumber
\\
&\quad& \pm \frac {\partial \psi}{\partial \nu} e^{-\eta_n\sqrt{\lambda}}. \label{normal derivatives on the boundary}
\end{eqnarray}

It follows from \eqref{assumption for induction} that for $x \in \overline{\mathcal N_-}$
\begin{equation}
\label{expansion with higher order curvatures ell}
\prod\limits_{j=1}^{N-1}\Bigl[1-\kappa_j(z(x))\delta(x)\Bigr] = 1 + \sum_{i=p}^{N-1} (-1)^i H_i(z(x)) (\delta(x))^i.
\end{equation}
We choose, for instance, $n = N-1$.   Let us show that for every $s \in \{ 0, \dots, p-2 \}$ as  $\delta(x) \to 0$
\begin{equation}
\label{estimates near the boundary}
\Delta A_s(x) = - 2^{-(s+1)} (-1)^p  (s+2)!\binom {p}{s+2} H_p(z(x)) \left(\delta(x)\right)^{p -2-s} + O\left(\left(\delta(x)\right)^{p -1-s}\right).
\end{equation}
By \eqref{expansion with higher order curvatures ell} and \eqref{A0}, we have that as  $\delta(x) \to 0$
$$
A_0(x) = 1 - \frac 12 (-1)^p H_p(z(x)) (\delta(x))^p +O\left(\left(\delta(x)\right)^{p +1}\right).
$$
Then, it follows from the first equality of \eqref{distance functions and principal curvatures} that as  $\delta(x) \to 0$
$$
\Delta A_0(x) = - \frac 12 (-1)^p H_p(z(x)) p(p-1)(\delta(x))^{p-2} + O\left(\left(\delta(x)\right)^{p -1}\right),
$$
which means that \eqref{estimates near the boundary} holds for $s = 0$. Suppose that \eqref{estimates near the boundary} holds for $s = q-1 \in \{ 0, \dots, p-2 \}$.
Then we have from \eqref{Aj} that
\begin{eqnarray*}
A_{q}(x) &=& \int_0^{\delta(x)}\left[- 2^{-(q+1)} (-1)^p  (q+1)!\binom {p}{q+1} H_p(z(x)) \tau^{p -1-q} + O\left(\left(\tau\right)^{p -q}\right)\right] \times
\\
&\quad&\exp\left(-\frac 12\,\int_\tau^{\delta(x)} \Delta \delta(x(\tau')) d\tau'\right)d\tau
\\
&=& - 2^{-(q+1)} (-1)^p  (q+1)!\binom {p}{q+1} H_p(z(x)) \int_0^{\delta(x)}  \tau^{p -1-q} d\tau + O\left(\left(\delta(x)\right)^{p -q+1}\right)
\\
&=& - 2^{-(q+1)} (-1)^p  q!\binom {p}{q} H_p(z(x)) {\delta(x)}^{p -q} + O\left(\left(\delta(x)\right)^{p -q+1}\right).
\end{eqnarray*}
Thus it follows from the first equality of \eqref{distance functions and principal curvatures} that as  $\delta(x) \to 0$
\begin{eqnarray*}
\Delta A_{q}(x) &=&  - 2^{-(q+1)} (-1)^p  q!\binom {p}{q} H_p(z(x)) (p-q)(p-q-1) {\delta(x)}^{p -q-2} + O\left(\left(\delta(x)\right)^{p -q-1}\right)
\\
&=& - 2^{-(q+1)} (-1)^p  (q+2)!\binom {p}{q+2} H_p(z(x)) {\delta(x)}^{p -q-2} + O\left(\left(\delta(x)\right)^{p -q-1}\right),
\end{eqnarray*}
which means that \eqref{estimates near the boundary} holds for $s = q$. Hence formula \eqref{estimates near the boundary} holds true for every $s\in \{0,\dots, p-2\}$.

Formula \eqref{estimates near the boundary}  implies that on $\partial\Omega$
$$
\Delta A_s = 0\quad \mbox{ for } s < p-2\ \mbox{ and }\ \Delta A_{p-2} = - 2^{-(p-1)} (-1)^p  p! H_p,
$$
and hence it follows from \eqref{key inequalities on the boundary more precise}  and \eqref{normal derivatives on the boundary} that on $\partial\Omega$ as $\lambda \to \infty$
\begin{equation}
\label{normal derivative from inside precise}
\sigma_s\frac {\partial w}{\partial\nu}\Big|_- = \frac {\sqrt{\sigma_s}\sqrt{\sigma_m}}{\sqrt{\sigma_s}+\sqrt{\sigma_m}}\left\{ \sqrt{\lambda} + p!2^{-p}(-1)^p(\sigma_s)^{\frac p2}H_p\lambda^{-\frac {p-1}2}\right\} + O\left( \lambda^{-\frac p2}\right).
\end{equation}

Next, as in the end of subsection \ref{subsec:vanishing mean curvature},  we proceed to consider $1-w$ on $\mathcal N_+$. By replacing $w,\ \sigma_s$ with $1-w,\ \sigma_m$, respectively and taking into account both the sign of $H_p$ and the normal direction to $\partial\Omega$,  by the same arguments we infer that on $\partial\Omega$ as $\lambda \to \infty$
\begin{equation}
\label{normal derivative from outside precise}
\sigma_m\frac {\partial w}{\partial\nu}\Big|_+= \frac {\sqrt{\sigma_s}\sqrt{\sigma_m}}{\sqrt{\sigma_s}+\sqrt{\sigma_m}}\left\{ \sqrt{\lambda} + p!2^{-p}(\sigma_m)^{\frac p2}H_p\lambda^{-\frac {p-1}2}\right\} + O\left( \lambda^{-\frac p2}\right).
\end{equation}
Here we used the fact that,  corresponding to the choice of the normal direction to $\partial\Omega$,  the sign of $H_p$ changes if $p$ is odd and it does not change if $p$ is even. Since $\sigma_s \not= \sigma_m$,  by combining \eqref{normal derivative from inside precise} and \eqref{normal derivative from outside precise} with the second equality of \eqref{interface condition for lambda} we conclude that on $\partial\Omega$
$$
H_p = O(1/\sqrt{\lambda})\ \mbox{ as } \lambda \to \infty,
$$
and hence $H_p$ must vanish  on $\partial\Omega$. Therefore we obtain that $H_i = 0$ on $\partial\Omega$ for every $i=1, \dots, N-1$. This means that all the principal curvaures of $\partial\Omega$ vanish and thus $\partial\Omega$ must be a hyperplane. 

Note that in this subsection \ref{subsec:vanishing principal curvatures} we used the fact that $\sigma_s \not=\sigma_m$.

%%%%%%%%%%%%%%%%%%%%%%%%%%%%%%%%%%%%%%%%%%%%%%%%%%%%%%
%%%%%%%%%%%%%%%%%%%%%%%%%%%%%%%%%%%%%%%%%%%%%%%%%%%%%%
%%%%%%%%%%%%%%             Section 4 begins        %%%%%%%%%%%%%%%%%%%%%%%%
%%%%%%%%%%%%%%%%%%%%%%%%%%%%%%%%%%%%%%%%%%%%%%%%%%%%%%
%%%%%%%%%%%%%%%%%%%%%%%%%%%%%%%%%%%%%%%%%%%%%%%%%%%%%%
\setcounter{equation}{0}
\setcounter{theorem}{0}

\section{Proof of Theorem \ref{th:characterization of a hyperplane and a helicoid in the heat equation}}
\label{section4}

Let $u$ be the solution of  problem \eqref{heat Cauchy one phase}.
From \eqref{stationary isothermic surface}  we see that $\partial\Omega$ is a stationary isothermic surface of $u$. Thus by \cite[Theorem 2.2, p.  4825]{MPStransaction2006} $\partial\Omega$ must be a real analytic hypersurface embedded in $\mathbb R^N$.  Hence Proposition \ref{modified formula on the boundary value} yields that  $k= \frac 12$.  Let $x \in \partial\Omega$. Then it follows from the explicit representation of $u$ via Gaussian kernel  that for every $t > 0$
\begin{eqnarray*}
\frac 12 &=& u(x, t) =  (4\pi t)^{-\frac N2}\int_{\mathbb R^N} {\mathcal X}_{\Omega^c}(\xi) e^{-\frac {|x-\xi|^2}{4t}}d\xi 
\\
&=& (4\pi t)^{-\frac N2}\int_0^\infty e^{-\frac {r^2}{4t}}\left(\int_{\partial B_r(x)} {\mathcal X}_{\Omega^c}(\xi) dS_\xi\right) dr
\\
&=& (4\pi t)^{-\frac N2}\int_0^\infty e^{-\frac {r^2}{4t}} | \Omega^c\cap\partial B_r(x) | dr,
\end{eqnarray*}
where $\Omega^c=\mathbb R^N \setminus \Omega$, $dS_\xi$ indicates the area element of the sphere $\partial B_r(x)$ and $| \Omega^c\cap\partial B_r(x) |$ does the $(N-1)$-dimensional Hausdorff measure of the set $\Omega^c\cap\partial B_r(x)$. Thus we infer that
$$
\int_0^\infty e^{-\frac {r^2}{4t}} \left(| \Omega^c\cap\partial B_r(x) | - \frac 12 |\partial B_r(x) |\right) dr =  0\ \mbox{ for every } t > 0.
$$
Since the Laplace transform is injective, we conclude that for each point $x \in \partial\Omega$ 
\begin{equation}
\label{uniformly dense Hausdorff measure version}
| \Omega^c\cap\partial B_r(x) | - \frac 12 |\partial B_r(x) | = 0\ \mbox{ for almost every } r > 0.
\end{equation}
Then the following formula also holds true:
\begin{equation}
\label{uniformly dense in volume}
\frac {|\Omega^c\cap B_r(x)|}{|B_r(x)|} = \frac 12\  \mbox{ for every } r > 0 \mbox{ and } x \in \partial\Omega,
\end{equation}
where the same symbol $| \cdot |$ indicates the $N$-dimensional Lebesgue measure of sets.

When $N \ge 2$, by \cite[Theorem 1.2, p.  4823]{MPStransaction2006} \eqref{uniformly dense in volume}
yields that $\partial\Omega$ must have zero mean curvature. Hence, when $N = 2$, $\partial\Omega$ must be a straight line, and when $N\ge 3$, $\partial\Omega$ must be a minimal hypersurface embedded in $\mathbb R^N$. 

In view of the sufficient regularity of $\partial\Omega$,  it follows from \eqref{uniformly dense Hausdorff measure version} that for every point  $p \in \partial\Omega$,  there exist  numbers $\delta_p > 0$ and $r_p > 0$ satisfying
\begin{equation}
\label{for Nitsche}
| \Omega^c\cap\partial B_r(x) | - \frac 12 |\partial B_r(x) | = 0\ \mbox{ for every } 0 < r < r_p \mbox{ and } x \in B_\delta(p)\cap \partial\Omega.
\end{equation}
When $N = 3$, by \cite[Theorem, p. 234]{Nanalysis1995},  \eqref{for Nitsche} yields that $\partial\Omega$ must be either a hyperplane  or a helicoid.   This completes the proof of Theorem \ref{th:characterization of a hyperplane and a helicoid in the heat equation}.  \qed

%%%%%%%%%%%%%%%%%%%%%%%%%%%%%%%%%%%%%%%%%%%%%%%%%%%%%%
%%%%%%%%%%%%%%%%%%%%%%%%%%%%%%%%%%%%%%%%%%%%%%%%%%%%%%
%%%%%%%%%%%%%%             Section 5 begins        %%%%%%%%%%%%%%%%%%%%%%%%
%%%%%%%%%%%%%%%%%%%%%%%%%%%%%%%%%%%%%%%%%%%%%%%%%%%%%%
%%%%%%%%%%%%%%%%%%%%%%%%%%%%%%%%%%%%%%%%%%%%%%%%%%%%%%

\setcounter{equation}{0}
\setcounter{theorem}{0}

\section{Appendices}
\label{section5}
First of all, let us give a proof of \eqref{just half on the helicoid}.  

 \subsection{Proof of \eqref{just half on the helicoid}}
\label{subsec: helicoid}
Let $\mathcal H \subset \mathbb R^3$ be the helicoid given by
$$
\left\{(x_1, x_2, x_3) =(\rho \cos s, \rho \sin s, s)\; :\; (\rho, s) \in \mathbb R^2\right\}.
$$
(See \cite[pp. 8--9]{CM2011} for the helicoid).   
Notice that $\mathcal H$ is the boundary of the following unbounded domain:
\begin{equation}\label{the one whose boundary is the helicoid}
\Omega=\left\{ (x_1,x_2,x_3)\in {\mathbb R}^3 \; :\;  x_2 \cos x_3 - x_1 \sin x_3>0 \right\}.
\end{equation}

We now introduce two isometries that are deeply related to the symmetries of $\mathcal H$. For $\alpha\in \mathbb R$ and $x=(x_1,x_2,x_3)\in {\mathbb R}^3$, we set: 
\begin{eqnarray}\nonumber
&k_\alpha(x)= \left( x_1  \cos\alpha  -  x_2\sin\alpha ,\ x_1\sin \alpha + x_2 \cos \alpha,\ x_3+\alpha  \right), \\
&g(x)= (x_1, -x_2, -x_3).
\end{eqnarray}
Here $k_\alpha$ is a \emph{screwing motion} obtained by rotation  of angle $\alpha$ in the $x_1$-$x_2$ plane, followed by a translation of length $\alpha$ in the $x_3$ direction. Notice that $\Omega$ and $\mathbb{R}^3\setminus \overline \Omega$ are preserved by the action of $k_\alpha$, while they get switched by that of $g$:
\begin{equation}\label{some symmetry properties} 
\begin{aligned}
k_\alpha(\Omega) = \Omega, \quad k_\alpha(\mathbb{R}^3\setminus \overline \Omega)=\mathbb{R}^3\setminus \overline \Omega,\\
g(\Omega)= \mathbb{R}^3\setminus \overline \Omega, \quad g(\mathbb{R}^3\setminus \overline \Omega)=\Omega.
\end{aligned}
\end{equation}
Finally, since  $x_2 \cos x_3 - x_1 \sin x_3=0$ for $x \in \mathcal{H}$, the restrictions of $g$ and $k_\alpha$ to $\mathcal H$ are related by the following formula:
\begin{equation}\label{sometimes they coincide}
g(x_1,x_2,x_3)=(x_1,-x_2, -x_3)=k_{-2 x_3}(x_1,x_2,x_3) \quad \text{ for all } (x_1,x_2,x_3)\in \mathcal{H}.
\end{equation}

Let $u =u(x,t)$ be the unique bounded solution of the following Cauchy problem for the heat diffusion equation:
\begin{equation}\label{heat Cauchy for helicoid}
  u_t = \Delta u\quad\mbox{ in }\  \mathbb R^3\times (0,+\infty) \ \mbox{ and }\ u\ ={\mathcal X}_{\mathbb R^3 \setminus \Omega}\ \mbox{ on } \mathbb R^3\times
\{0\},
\end{equation}
where $\Omega$ is the unbounded domain defined in \eqref{the one whose boundary is the helicoid}.
Moreover, for arbitray real $\alpha$, define the following functions:
$$
v_\alpha(x,t)=u\left( k_\alpha(x), t \right)\quad \text{and}\quad w(x,t)=u\left( g(x), t  \right)\quad \text{for} \quad (x,t)\in\mathbb R^3\times(0,\infty).
$$
Since both $k_\alpha$ and $g$ are isometries, by \eqref{some symmetry properties} we deduce that $v_\alpha$ and $w$ 
are bounded solutions of the following Cauchy problems.
\begin{eqnarray}\label{v_al eq}
 (v_\alpha)_t = \Delta v_\alpha\quad\mbox{ in }\  \mathbb R^3\times (0,+\infty) \ &\mbox{ and }\ & v_\alpha\ ={\mathcal X}_{\mathbb R^3 \setminus \Omega}\ \mbox{ on } \mathbb R^3\times\{0\},\\
\label{w eq} w_t = \Delta w\quad\mbox{ in }\  \mathbb R^3\times (0,+\infty) \ &\mbox{ and }\ & w\ ={\mathcal X}_{\Omega}\ \mbox{ on } \mathbb R^3\times\{0\}.
\end{eqnarray}
In particular, unique solvability of the Cauchy problems above yields 
\begin{equation}\label{symmetry by uniqueness}
v_\alpha=u \quad \text{and}\quad u+w=1\quad \text{in }\mathbb R^3\times(0,\infty),\quad \text{for all }\alpha\in\mathbb R .
\end{equation}

Fix now an arbitrary pair $(x,t)\in \mathcal H \times (0,\infty)$ and choose $\alpha=-2x_3$. By combining both identities in \eqref{symmetry by uniqueness} with \eqref{sometimes they coincide} we get the following chain of equalities.
$$
1=u(x,t)+u(g(x),t)= u(x,t)+u\left(k_{-2x_3}(x), t\right)=2u(x,t).
$$
That is, $u(x,t)=1/2$ for all $(x,t)\in \mathcal H \times (0,\infty)$. We have therefore proved \eqref{just half on the helicoid} when $N=3$. The case $N\ge 4$ follows by separation of variables. \qed

\subsection{A maximum principle for unbounded domains}
\label{subsec: max pple}

For convenience, we quote  a maximum principle together with its proof for an elliptic equation in unbounded domains in $\mathbb R^N$ from \cite[Proposition A.3]{SarXiv2019}.
%%%%%%%%%%%%%%%%%%%%%%%%%%%%%%%%%%%%%%%%%
%%%%%%%%%%%%%%%%%%%%%%%%%%%%%%%%%%%%%%%%%
%%%%%     Proposition A.1 begins       %%%%%%%%
%%%%%%%%%%%%%%%%%%%%%%%%%%%%%%%%%%%%%%%%%
%%%%%%%%%%%%%%%%%%%%%%%%%%%%%%%%%%%%%%%%%

\begin{proposition}
\label{prop:maximum principle on unbounded domains}
Let $D \subset \mathbb R^N$ be an unbounded domain, and let  $\sigma=\sigma(x)\ (x\in D)$ be 
a general conductivity satisfying
$$
0 < \mu \le \sigma(x) \le M\ \mbox{ for every } x \in \mathbb R^N,
$$
where $\mu, M$ are positive constants.  Assume that $w \in H^1_{loc}(D)\cap L^\infty(D)\cap C^0(\overline{D})$ satisfies
$$
-\mbox{\rm div}(\sigma \nabla w) + \lambda w \ge 0\ \mbox{ in } D\ \mbox{ and }\ w \ge 0\ \mbox{ on } \partial D
$$
for some constant $\lambda > 0$. Then $w \ge 0$ in $D$,  and moreover, either $w > 0$ in $D$ or $w \equiv 0$ in $D$.
\end{proposition}

%%%%%%%%%%%%%%%%%%%%%%%%%%%%%%%%%%%%%%%%%
%%%%%%%%%%%%%%%%%%%%%%%%%%%%%%%%%%%%%%%%%
%%%%%     Remark A.2 begins       %%%%%%%%
%%%%%%%%%%%%%%%%%%%%%%%%%%%%%%%%%%%%%%%%%
%%%%%%%%%%%%%%%%%%%%%%%%%%%%%%%%%%%%%%%%%
\begin{remark}
When $D$ is bounded, this proposition is well known and holds true  for every $\lambda \ge 0$. However, when $D$ is unbounded, this proposition is not true for $\lambda = 0$. Indeed, a  counterexample is given in {\rm \cite[p. 37]{ABR2001Springer}}, where $N \ge 3,\ D = \{ x \in \mathbb R^N : |x| > 1 \},\ \sigma(x)\equiv 1$ and  $w(x) = |x|^{2-N}-1$.
\end{remark}

\noindent
{\it Proof of Proposition \ref{prop:maximum principle on unbounded domains}.} 
Define $v = v(x)$ by
$$
v(x) = e^{-\delta|x|} w(x)\ \mbox{ for } x \in \overline{D},
$$
where  $\delta > 0$ is a constant which will be chosen later.
Then $v \in H^1_{loc}(D)\cap L^\infty(D)\cap C^0(\overline{D})$ and moreover 
\begin{equation}
\label{limit zero at infinity}
\lim_{|x| \to \infty} v(x) = 0,
\end{equation}
since $w \in L^\infty(D)$.
For every $ \varepsilon > 0$,  we consider a nonnegative function
$$
\varphi(x) = \max\{ -\varepsilon - v(x), 0 \}\  \mbox{ for } x \in \overline{D}.
$$
Since $v \in H^1_{loc}(D)\cap L^\infty(D)\cap C^0(\overline{D})$ and $v \ge 0\ \mbox{ on } \partial D$, it follows from \eqref{limit zero at infinity} that $\varphi$ is compactly supported in $D$ and $\varphi \in H^1_0(D)$, and hence $e^{-2\delta|\cdot|}\varphi(\cdot)  \in   H^1_0(D)$. Therefore we obtain
\begin{eqnarray}
 0 &\le& \int\limits_D\left\{\sigma(x) \nabla w(x) \!\cdot \!\nabla\!\left(\varphi(x) e^{-2\delta|x|}\right) + \lambda w(x) \varphi(x)e^{-2\delta|x|}\right\} dx \nonumber
\\
&=& \int\limits_{D\cap\{v<-\varepsilon\}}\!\!\!\!\!\!\! \sigma e^{-\delta|x|}\left\{\left(\delta{v}\frac x{|x|} +\nabla{v}\right)\!\cdot\!\left(\nabla\varphi-2\delta\varphi\frac x{|x|}\right) + \frac\lambda\sigma{ v}\varphi \right\} dx. \label{the integral including varepsilon and delta}
\end{eqnarray}
Notice  that 
$$
\varphi(x) = 
\begin{cases}
-\varepsilon - v(x)   &\mbox{if }\ v(x) < -\varepsilon, \\
0  &\mbox{if }\ v(x) \ge -\varepsilon,
\end{cases} 
\quad\mbox{ and }\quad
\nabla \varphi(x) = 
 \begin{cases}
-\nabla v(x)   &\mbox{if }\ v(x) < -\varepsilon, \\
0  &\mbox{if }\ v(x) \ge -\varepsilon.
\end{cases} 
$$
By setting
$$
I = \sigma^{-1} e^{\delta |x|} \times \mbox{ the integrand of the integral \eqref{the integral including varepsilon and delta}},
$$
we have
\begin{eqnarray*}
I &=& -|\nabla v|^2-\frac\lambda\sigma v^2+2\delta^2v^2+\delta v\frac x{|x|} \cdot\nabla v + \varepsilon\left(2\delta^2 v+2\delta\frac x{|x|} \cdot\nabla v-\frac\lambda\sigma v\right)
\\
&\le & -\left\{ 1-\delta\left(\frac 12+\varepsilon\right)\right\}|\nabla v|^2 - \left\{ \frac\lambda\sigma\left(1-\frac\varepsilon 2\right) -\left(2\delta^2+\frac \delta 2\right) \right\} v^2 + \varepsilon\left(\frac \lambda {2\sigma}+\delta\right).
\end{eqnarray*}
Here we have used Cauchy's inequality $2ab \le a^2+b^2$ and the fact that $v< 0$ in the integrand of \eqref{the integral including varepsilon and delta}.  Therefore, since $0 < \mu \le \sigma(x)\le M$, we can choose $\delta > 0 $ sufficiently small to obtain that if $0 < \varepsilon < 1$ then 
$$
I \le -\frac 14\left( |\nabla v|^2 + \frac\lambda M v^2\right) + \varepsilon\left(\frac \lambda {2\mu}+\delta\right)
$$
and hence
$$
\mu\!\!\!\!\!\!\!\!\! \int\limits_{D\cap\{v<-\varepsilon\}}\!\!\!\!\!\!\! e^{-\delta|x|}\left( |\nabla v|^2 + \frac \lambda M v^2\right) dx \le M\varepsilon\left(\frac {2\lambda}\mu+4\delta\right)\int\limits_{D} e^{-\delta |x|} dx.
$$
By choosing a sequence $\{\varepsilon_n\}$ with $\varepsilon_n\downarrow 0$ as $n \to \infty$ and letting $n \to \infty$, we conclude that
$$
\int\limits_{D\cap\{v<0\}}\!\!\!\!\! e^{-\delta|x|}\left( |\nabla v|^2 + \frac\lambda M v^2\right) dx = 0
$$
and hence $v \ge 0$ in $D$.  Therefore $w \ge 0$ in $D$. Once this is shown, the last part follows from the strong maximum principle (see \cite[Theorem 8.19, pp. 198--199]{GT1983}). \qed

\end{document}